\documentclass[a4paper, 10pt,russian]{article}
\usepackage[cp866]{inputenc}
\usepackage[russian]{babel}
\usepackage{color}
\usepackage{ifthen}
\usepackage{amsfonts,amsmath,amssymb,theorem}
\usepackage[unicode]{hyperref} 

\makeatletter
\def\@seccntformat#1{\csname the#1\endcsname.\ } 
\def\@biblabel#1{#1.} 
\makeatother

\date{}

\newcommand\proofr{\par{\it Доказательство.\,~}}

\def\proofend{$\blacktriangle$\vspace{0.3em}\par}

{\par\addvspace{1mm}{\it Proof\hspace{1.0ex}{#1}.} }%
{\quad$\blacktriangle$\par\addvspace{1mm}}

    \newif\ifNoRemark
    \def\addtheorem#1#2#3#4{ 
    \ifthenelse{\expandafter\isundefined\csname the#2\endcsname}{\newcounter{#2}}{}
    \newenvironment{#1}[1][\global\NoRemarktrue]
     {\par\addvspace{2mm}\noindent  
       \refstepcounter{#2}{\bf #3~\csname the#2\endcsname
      \vphantom{##1}\ifNoRemark.\ \else\ (##1).\fi}\begingroup #4}%
     {\endgroup\par\addvspace{1mm}\global\NoRemarkfalse}
    \expandafter\newcommand\csname b#1\endcsname{\begin{#1}}
    \expandafter\newcommand\csname e#1\endcsname{\end{#1}}
    }

\addtheorem{pro}{pro}{Утверждение}{\sl}
\addtheorem{theorem}{theorem}{Теорема}{\sl}
\addtheorem{lemma}{lemma}{Лемма}{\sl}
\addtheorem{corol}{corol}{Следствие}{\sl}
\addtheorem{remark}{rmrk}{Замечание}{\sl}

\newcommand{\supp}{{\rm supp}}
\newcommand{\rank}{{\rm rank}}

\title{О спектре мощностей\\ и числе латинских битрейдов порядка~$3$%
\thanks{Работа выполнена за счет грантов Российского научного фонда 14-11-00555  (разделы 2, 3, 4.1, 4.5, 5.2) и
18-11-00136 (разделы 4.2--4.4, 4.6, 5.1). }\\
(On the cardinality spectrum\\ and the number of latin bitrades of order~$3$)%
\thanks{The results of this work were presented in parts at the  International Conference and PhD-Master Summer School
``Groups and Graphs, Metrics and Manifolds'' G2M2,
Yekaterinburg, July 22--30, 2017, and the International Workshop on Algebraic Combinatorics, Hefei, China, November 22--25, 2018.}
}
\author{D. S. Krotov, V. N. Potapov\\
\emph{Sobolev Institute of Mathematics}
\\krotov@math.nsc.ru; vpotapov@math.nsc.ru}

\begin{document}

\maketitle

\begin{abstract}
By a (latin) unitrade, we call a set of vertices of the Hamming graph 
that is intersects with every maximal clique in $0$ or $2$ vertices.
A bitrade is a bipartite unitrade, that is, a unitrade splittable into two independent sets.
We study the cardinality spectrum of the bitrades in the Hamming graph $H(n,k)$ with $k=3$ (ternary hypercube)
and the growth of the number of such bitrades as $n$ grows.
In particular, we determine all possible (up to  $2.5\cdot 2^n$) and large (from $14\cdot 3^{n-3}$)
cardinatities of bitrades and prove that the cardinality of a bitrade is compartible to $0$ or $2^n$ modulo $3$
(this result has a treatment in terms of a ternary code of Reed--Muller type). A part of the results is valid for any $k$.
We prove that the number of nonequivalent bitrades is not less than
$2^{(2/3-o(1))n}$ and is not greater than $2^{\alpha^n}$, $\alpha<2$, as $n\rightarrow\infty$.

Унитрейдом (латинским) называется подмножество вершин графа Хэмминга, которое либо
пересекается по двум вершинам, либо совсем не пересекается с любой
максимальной кликой. 
Битрейдом называется двудольный, т.\,е.
разделяющийся на два независимых подмножества, унитрейд. 
В статье
исследуется спектр мощностей битрейдов в графе Хэмминга $H(n,k)$ при $k=3$
(троичном гиперкубе) и рост числа таких битрейдов с ростом $n$. В частности, определены все возможные  малые
(до $2.5\cdot 2^n$) и большие (от $14\cdot 3^{n-3}$) мощности битрейдов размерности $n$ и доказано,
что мощность битрейда принимает значения только сравнимые с $0$ или
$2^n$ по модулю $3$ (этот результат имеет трактовку в терминах троичного кода типа Рида--Маллера). Часть результатов применима для произвольного $k$. Доказано, что число неэквивалентных битрейдов
не меньше $2^{(2/3-o(1))n}$ и не больше $2^{\alpha^n}$, $\alpha<2$,
при $n\rightarrow\infty$.
\end{abstract}


\section{Введение}\label{s:intro}

Для комбинаторных объектов (конфигураций) различного типа
оказывается полезно  вводить понятие битрейда, так чтобы определение
битрейда не опиралось непосредственно на определение исходных
объектов, но  включало всевозможные разности (например,
симметрические) объектов этого типа (см. \cite{Kro16}).
 Битрейды отражают возможную разницу между двумя комбинаторными
конфигурациями одного и того же типа, что важно при перечислении,
описании и исследовании свойств комбинаторных конфигураций. Известны
исследования битрейдов (или трейдов) комбинаторных  блок-дизайнов
\cite{HK,KMT,Kro18,KMP:16:trades,GKKK},  латинских квадратов \cite{Cav},
частично упорядоченных множеств \cite{Cho},  совершенных кодов
\cite{AS,Ost:2012:switching}, корреляционно-иммунных и бент- функций \cite{PP12}. В
настоящей работе рассматриваются латинские битрейды, соответствующие
МДР-кодам с расстоянием $2$, или (что эквивалентно) латинским
гиперкубам, или полиадическим квазигруппам. 
В работе исследован спектр возможных мощностей
битрейдов и получены оценки их числа.

Перейдём к формальным определениям. Пусть $Q_k=\{0,\dots,k-1\}$.
Определим расстояние Хэмминга $d(u,v)$ как число несовпадающих
компонент в наборах $u,v\in Q_k^n$. Метрическое пространство
$(Q^n_k,d)$, а также граф $\Gamma Q_k^n$ расстояний $1$ на множестве
вершин $Q_k^n$, называется $k$-ичным $n$-мерным гиперкубом или
графом Хэмминга. Весом вершины $u\in Q_k^n$ называется
$\mathrm{wt}(u)=d(u,\overline{0})$, 
где $\overline{0}$ (далее также используются обозначения $\overline{1}$, $\overline{-1}$) --- 
набор из $n$ нулей (соответственно, $1$-ц, $-1$-ц). 
Гранью в $Q_k^n$ называется
 множество вершин гиперкуба, полученное фиксацией значений одной или
нескольких координат. Множество $U\subset Q_k^n$ называется {\it
унитрейдом (размерности $n$)}, если мощности его пересечений с
одномерными гранями (максимальными кликами в $\Gamma Q_k^n$)
принимают только два значения $0$ и $2$. Обычно битрейдом называется
пара $\{U_0,U_1\}$, состоящая из двух независимых долей двудольного
унитрейда $U=U_0\cup U_1$. Но нам будет удобно называть {\it
битрейдом} такой унитрейд $U\subset Q_k^n$, что
  подграф
$\Gamma U$, порождённый множеством вершин $U$, является двудольным.
В двумерном случае ($n=2$) из теоремы Кёнига следует, что любой
унитрейд является битрейдом. При $n\geq 3$ и $k\geq 3$ имеются
унитрейды, которые не являются битрейдами.

Рассмотрим соответствие между данным выше определением и общим
понятием битрейда. {\it МДР-кодом} называется подмножество гиперкуба
$Q_k^n$, пересекающееся с каждой гранью фиксированной размерности
$r$ ровно по одному элементу. Нетрудно видеть что МДР-коды --- это
коды с минимальным расстоянием между вершинами $r+1$ и максимальной
для этого кодового расстояния мощностью $k^{n-r}$, т.\,е. коды
достигающие границы Синглтона. В данном контексте нас интересуют 
только МДР-коды с кодовым расстоянием $2$, т.\,е. когда $r=1$. (Функция, выражающая значение одной из координат вершин такого кода через значения $n-1$ оставшихся называется латинским $(n-1)$-кубом, в случае $n=3$ латинским квадратом, а алгебраическая система с такой функцией в качестве операции --- $(n-1)$-арной квазигруппой). Из
определений видно, что симметрическая разность двух МДР-кодов
является битрейдом.


Группа изометрий гиперкуба $Q_k^n$ порождается группой перестановок
координат и группой изотопий, т.\,е. перестановок элементов $Q_k$ в
каждой координате. В случае $k=3$ группа изометрий $Q_3^n$
состоит только из аффинных преобразований гиперкуба $Q_3^n$ как
$n$-мерного векторного пространства над полем $\mathrm{GF}(3)$. Подмножества
гиперкуба, которые можно перевести друг в друга изометриями
пространства называются {\it эквивалентными}. {\it Ретрактом }
множества $U\subset Q_k^n$ будем называть подмножество гиперкуба
$Q_{k}^{n-1}$, полученное как пересечение множества $U$ с некоторой
гранью размерности $n-1$.
 Из определений следуют

\bpro\label{probit00}
 Любой ретракт унитрейда (битрейда, {\rm МДР}-кода)
является унитрейдом (битрейдом, {\rm МДР}-кодом) в гиперкубе меньшей
размерности.
 \epro

\bpro\label{probit000}
 Образ унитрейда (битрейда, {\rm МДР}-кода) при изометрии гиперкуба
является унитрейдом (битрейдом, {\rm МДР}-кодом).
 \epro

В статье мы почти полностью ограничиваемся рассмотрением битрейдов в
$Q^n_3$. Этот случай представляется нам  ключевым, поскольку любой
троичный битрейд можно изометрично вложить в гиперкуб большего
порядка $Q^n_k$, $k\geq 4$. Причём, если в гиперкубах  порядка $k$,
$k\geq 4$, имеются унитрейды и битрейды, состоящие из нескольких
компонент связности, то в троичном случае никакой унитрейд не может
включать другой непустой унитрейд как подмножество и тем более
никакой троичный унитрейд нельзя представить в виде объединения
непересекающихся унитрейдов. Для исследования унитрейдов в статье
используются методы линейной алгебры (\S\,\ref{s:lin}, \S\,\ref{ss:ter-lin}), теории булевых
функций (\S\,\ref{s:bool}) и теории кодирования (\S\,\ref{s:bool}, \S\,\ref{ss:rm}). В частности,
показана связь задачи описания троичных битрейдов с задачей
нахождения полиномиальной сложности булевой функции \cite{RESM},
\cite{VK}. Известна также связь троичных битрейдов с почти
уравновешенными булевыми функциями (утверждение \ref{proPot12}).

В разделе \S\,\ref{s:spectr} исследуется спектр мощностей троичных битрейдов.
Показано, что мощность любого битрейда размерности $n$ сравнима с
$0$ или $2^n$ по модулю три. Минимальная мощность непустого битрейда
(не только троичного) размерности $n$ равна $2^n$. Все возможные
мощности битрейдов, не превышающие $2\cdot 2^n$, были известны ранее
(см. \cite{Pot13}). В настоящей работе указана связь  спектра
мощностей троичных битрейдов с весовым спектром двоичного кода
Рида--Малера. А также определены все возможные мощности битрейдов
размерности $n$ до $2.5\cdot 2^n$ (теорема~\ref{th:2.5}). Битрейдом
максимальной мощности $2\cdot 3^{n-1}$ является пара
непересекающихся МДР-кодов. Как известно (см., например,
\cite{KP09}) имеется единственный (с точностью до эквивалентности)
битрейд такой мощности. Отметим, что уже для порядка $4$  имеется
дважды экспоненциальное число неэквивалентных битрейдов максимальной
мощности $2\cdot 4^{n-1}$ (см. \cite{KP09}, \cite{KP11}).

Одной из основных задач исследования битрейдов является определение
их числа как функции от размерности $n$ и порядка $k$. Благодаря
тому, что битрейды соответствуют разностями комбинаторных объектов,
исследование разнообразия битрейдов открывает перспективы
исследования разнообразия исходных объектов и оценки их числа.
 Для
латинских гиперкубов  (порядка больше $4$) размерности $n$, с
которыми связаны исследуемые нами битрейды, до сих пор неизвестен
даже порядок роста логарифма числа при $n\rightarrow\infty$ (см.
\cite{KP11}). В \cite{Pot13} была получена почти экспоненциальная
$(e^{\Omega(\sqrt n)})$ асимптотическая нижняя оценка числа
неэквивалентных битрейдов в $Q_3^n$. В \S\,\ref{ss:lb} доказана нижняя оценка
$2^{(2/3-o(1))n}$  числа неэквивалентных битрейдов в $Q_3^n$ при
$n\rightarrow\infty$ (теорема~\ref{th:lb}). В \S\,\ref{ss:ub} получена верхняя оценка
вида $2^{\alpha^n}$, $\alpha<2$, для числа битрейдов в $Q_3^n$
(теорема \ref{corlast}), которая существенно улучшает тривиальную
верхнюю оценку $2^{2^n}$. Однако, вопрос о скорости роста числа
троичных битрейдов размерности $n$ остаётся открытым в самом сильном
смысле: является ли эта функция экспоненциальной или дважды
экспоненциальной по $n$? Результаты численного эксперимента  для
малых $n$, представленые в таблице в разделе~\ref{ss:tables}, 
показывают рост быстрее экспоненцияльного, 
но вывод о линейном росте двойного логарифма 
числа представляется пока поспешным.

Метод доказательства верхней оценки связан с тем, что в гиперкубе
$Q^n_3$ при $n\geq 7$ удалось указать множество мощности cтрого
больше $3^n-2^n$, которое не включает в себя ни одного битрейда и
даже  симметрических разностей двух битрейдов. Для более известной
задачи о максимальном подмножестве гиперкуба без арифметической
прогрессии (что в троичном гиперкубе совпадает с подмножеством,
не включающем аффинного $1$-мерного подпространства) недавно была
получена асимптотическая верхняя оценка вида $o(\alpha^n)$, где
$\alpha<3$ (см. \cite{EG}). Нахождение в троичном гиперкубе
подмножества максимальной мощности, которое не включает битрейдов,
остаётся актуальной задачей.

Перечисление битрейдов в $Q_3^7$ проводилось 
на кластере Информационно-вычислительного центра
Новосибирского государственного университета (ИВЦ НГУ).


\section{Линейные пространства}\label{s:lin}

Пусть $\mathbb{F}$ --- некоторое поле. Рассмотрим множество функций
$\{g:Q_k^n\rightarrow \mathbb{F}\}$ как векторное пространство над
полем $\mathbb{F}$. Обозначим через $\mathbb{V}_{n,k}(\mathbb{F})$
подпространство, состоящее из функций, сумма значений которых по
любой $1$-мерной грани (максимальной клике в графе $\Gamma Q_k^n$)
равна $0$. Рассмотрим битрейд $B\subset Q_k^n$. Ему соответствует
функция $b[B]:Q_k^n\rightarrow \mathbb{F}$, которая на одной доле
битрейда принимает значение $1$, на другой --- значение $-1$ (для поля
характеристики $2$ имеем $-1=1$) и в остальных вершинах равна $0$. Очевидно,
$b[B]\in \mathbb{V}_{n,k}(\mathbb{F})$. Характеристическая функция
унитрейда содержится в $\mathbb{V}_{n,k}(\mathbb{F})$ в случае,
когда поле $\mathbb{F}$ имеет характеристику $2$.

 Введём следующий частичный порядок $\preceq$
на $Q_k$: элемент $k-1$ будем считать максимальным, а все остальные
элементы из $Q_k$ несравнимыми друг с другом. Распространим
частичный порядок на $Q_k^n$. Пусть $(x_1,\dots,x_n),
(y_1,\dots,y_n)\in Q_k^n$.
   Введём обозначение $(x_1,\dots,x_n)\preceq (y_1,\dots,y_n)$, если для любого
    $i\in \{1,\dots,n\}$ верно, что $x_i\preceq y_i$.
Отметим, что множество $G_y=\{x\in Q_{k}^n\ |\ x \preceq y\}$ является
гранью
 гиперкуба $Q_{k}^n$ размерности  равной числу
 символов $k-1$ в наборе $y$.

Покажем, что $\dim \mathbb{V}_{n,k}(\mathbb{F})=(k-1)^n$. Пусть
$f\in\mathbb{V}_{n,k}(\mathbb{F})$. Сумма значений функции $f$ по
каждой грани любой ненулевой размерности равна нулю, поэтому
\begin{equation}\label{eqbit1}
f(y)=-\sum\limits_{x\preceq y}f(x).
\end{equation}
Следовательно, для определения функции
$f\in\mathbb{V}_{n,k}(\mathbb{F})$ необходимо и достаточно
определить её на всех минимальных элементах, т.\,е. на $Q_{k-1}^n$.
Построим семейство линейно независимых функций той же мощности.
Пусть $x\in Q_{k-1}^n$. Рассмотрим множество $B_x=\{y\in Q_{k}^n \
|\ x\preceq y\}$. Нетрудно видеть, что граф $\Gamma B_x$
изоморфен булеву гиперкубу $\Gamma Q_2^n$, в частности $B_x$
является битрейдом. Здесь набору $z\in Q_2^n$ соответствует вершина
$y\in B_x$, у которой равны $k-1$ координаты, соответствующие
единицам в наборе $z$,  а координаты, соответствующие нулям в наборе
$z$, такие же как в наборе $x$. Таким образом верно равенство
$\mathrm{wt}(y-(k-1)\overline{1})=n-\mathrm{wt}(z)$. Соответствующую этому битрейду
функцию можно задать явной формулой
$b_x(y)=(-1)^{\mathrm{wt}(y-(k-1)\overline{1})}\chi_{ _{B_x}}(y)$. Поскольку
$\supp(b_x)\cap Q_{k-1}^n=\{x\}$, набор функций $\{b_x \ |\ x\in
Q_{k-1}^n\}$ является базисом в $\mathbb{V}_{n,k}(\mathbb{F})$.

Следующее утверждение нетрудно доказать по индукции (см.
\cite{Pot13})

\bpro\label{probit06} Пусть $f\in \mathbb{V}_{n,k}(\mathbb{F})$ и
$\supp(f)\neq \varnothing$. Тогда (a) $|\supp(f)|\geq 2^n$;\\ (b) если
$|\supp(f)|= 2^n$, то граф $\Gamma(\supp(f))$ изоморфен булеву
гиперкубу $\Gamma Q_2^n$.
  \epro

Нетрудно видеть, что имеется всего ${k \choose 2}^n$ вариантов
выбрать унитрейд (битрейд) с носителем мощности $2^n$.
 Как показано выше, базис  пространства
$\mathbb{V}_{n,k}(\mathrm{GF}(2))$ можно составить из характеристических
функций булева гиперкуба (точнее множеств, индуцирующих подграф, изоморфный булевому гиперкубу размерности $n$).
Поскольку при $k>2$ число таких множеств больше размерности пространства,
унитрейды не единственным образом представляются в виде линейной
комбинации над $\mathrm{GF}(2)$ булевых гиперкубов. Минимальное число булевых
гиперкубов в таком представлении унитрейда $U$ будем называть {\it
рангом } унитрейда и обозначать $\rank(U)$.

Рассмотрим подробнее троичный гиперкуб.  Любой элемент пространства
$\mathbb{V}_{n,3}({\mathrm{GF}(2)})$ является унитрейдом, поскольку чётное
число единиц из трёх возможных в $1$-мерной грани равняется $0$ или $2$.
Размерность пространства $\mathbb{V}_{n,3}({\mathrm{GF}(2)})$ равна $2^n$,
поэтому в гиперкубе $Q_3^n$ имеется ровно $2^{2^n}$ различных
унитрейдов. Методом индукции по размерности $n$ нетрудно доказать

\bpro\label{probit05} (a) Любая пара непустых унитрейдов в $Q^n_3$
имеет непустое
пересечение.\\
(b) Если в $Q^n_3$  непустой унитрейд  является подмножеством другого унитрейда,
то эти унитрейды совпадают.\\
(c)  Граф $\Gamma U$ унитрейда в $Q^n_3$ связен. \epro

Из (c)  следует, что если унитрейд $U\subset Q^n_3$ является
битрейдом, то он разделяется на две доли однозначно.


\section{Булевы функции}\label{s:bool}

Пусть $f:Q^n_2\rightarrow Q_2$ --- некоторая булева функция.
  Определим функцию $U[f]:Q_3^n\rightarrow \{0,1\}$
  равенством
$
 U[f](y)=
  \bigoplus\limits_{x\preceq y}f(x) $ (здесь и далее операция $\oplus$ обозначает сложение в двоичном поле). Из определения видно, что $U[f]|_{Q_2^n}=f$.
   Из совпадения определения функции $U[f]$ и
   формулы (\ref{eqbit1}) над полем $\mathrm{GF}(2)$ следует,
  что $U[f]$ --- характеристическая функция унитрейда в $Q_3^n$.
  Более того, из вышеизложенного следует, что между булевыми функциями
  и троичными унитрейдами имеется взаимно однозначное соответствие.

Пусть $x\in Q_2^n$. Введём обозначения $x_i^1=x_i$,
$x_i^{-1}=x_i\oplus 1$, $x_i^0=1$ и если $x=(x_1,\dots,x_n)$, то
$x^v=x_1^{v_1}\dots x_n^{v_n}$, где $v\in Q_3^n= \{0,\pm 1\}^n$.

{\it Полиномиальным представлением} булевой функции $f$ называется
формула вида\\ $f(x)=f^A(x_1,\dots,x_n)=\bigoplus\limits_{v\in
A\subset Q^n_3}x^v$. Минимальное количество слагаемых в этом
представлении ($\min|A|$) называется {\it полиномиальной сложностью}
функции $f$ (см. \cite{VK}).

Обозначим $\{0,\pm 1\}_0=\{0,1\}$, $\{0,\pm 1\}_1=\{1,-1\}$,
$\{0,\pm 1\}_{-1}=\{0,-1\}$ и $\{0,\pm 1\}_v=\{0,\pm
1\}_{v_1}\times\cdots\times\{0,\pm 1\}_{v_n}$. Вложенные в  $Q^n_3$
булевы гиперкубы исчерпываюся кубами вида $\{0,\pm 1\}_v$.  Нетрудно
видеть, что сужение характеристической функции гиперкуба $\{0,\pm
1\}_v$ на булев подкуб $\{0,1\}^n$ совпадает с мономом $x^v$, т.\,е.
$\chi_{_{\{0,\pm 1\}_v}}(x)=x^v$ при $x\in Q^n_2$. Тогда
$U[f^A]=\bigoplus\limits_{v\in A\subset Q^n_3}\chi_{_{\{0,\pm
1\}_v}}$ и  $\rank(U[f])$ равен полиномиальной сложности функции
$f$. Проблема нахождения  представления булевой функции минимальной
полиномиальной сложности (minimization of Exclusive or Sum of
Products) рассматривается, например, в \cite{RESM}. Известно
(см.\cite{VK}), что для размерности $5$ максимальная сложность
булевой функции равна $9$, для размерности $6$ --- равна $15$ и существуют
 булевы функции от $7$ аргументов с полиномиальной сложностью $24$.

Из определения полиномиальной сложности следует, что полиномиальная
сложность булевой функции не больше суммы сложностей двух её
подфункций, полученных фиксацией $0$ и $1$ некоторой переменной.
Отсюда имеем

\begin{corol}\label{corrang}
Ранг унитрейда не превосходит суммы рангов двух его различных
ретрактов по произвольной координате.
\end{corol}

Полиномиальное представление булевой функции неоднозначно. Но если
не использовать один из  операторов  $x^0$, $x^1$ или $x^{-1}$, то
представление приобретает однозначность. В предыдущем разделе
рассматривался базис в пространстве $f\in \mathbb{V}_{n,3}(\mathrm{GF}(2))$,
соответствующий операторам $x^1$ и $x^{-1}$. Если исключить оператор
$x^{-1}$ ("отрицание") мы приходим к базису из сложения и умножения
над $\mathrm{GF}(2)$. А именно, каждая булева функция $f:Q_2^n\rightarrow
Q_2$ может быть единственным образом представлена в виде {\it
многочлена Жегалкина} (в {\it алгебраической нормальной форме})
$$f(x_1,\dots,x_n)=\bigoplus\limits_{y\in
Q^n_2}G[f](y)x_1^{y_1}\dots x_n^{y_n}, $$ где $G[f]:Q_2^n\rightarrow
Q_2$ --- булева функция.

{\it Алгебраической степенью } булевой функции $f$ называется
максимальная степень слагаемого в её многочлене Жегалкина, т.\,е.
${\rm deg}\,f=\max\limits_{G[f](y)=1} \mathrm{wt}(y)$.


Справедливо следующее

\bpro\label{stsp0} Для любой булевой функции $f$ верно равенство
$G[f](y)=\bigoplus\limits_{x\in Q_2^n, x\preceq y}f(x)$. \epro

Таким образом, $G[f]$ является преобразованием Мёбиуса функции $f$
над полем $\mathrm{GF}(2)$. Поскольку $f(x)=\bigoplus\limits_{y\in Q_2^n,
x\preceq y}G[f](y)$, имеем равенство $G[G[f]]=f$ для любой булевой
функции $f$. Из определений преобразования Мёбиуса и оператора
$U[\cdot]$ видно, что $U[f]|_{\{0,2\}^n}$  является преобразованием
Мёбиуса булевой функции $f$.

Из утверждения  \ref{stsp0} непосредственно следует известное

\bpro\label{stsp1} Булева функция $f:Q_2^n\rightarrow Q_2$ имеет
чётное число единиц во всех гранях  размерности не меньше $m$ тогда
и только тогда, когда ${\rm deg}\,f\leq m-1$. \epro

Наборы значений булевых функций $f:Q_2^n\rightarrow Q_2$ можно
рассматривать  как элементы булева куба размерности $2^n$. Множество
наборов значений булевых функции алгебраической степени не выше $m$
называется кодом Рида--Маллера $\mathcal{R}(m,n)$ в
$Q_2^{2^n}$. Известно, что минимальный  вес ненулевого кодового
набора $\mathcal{R}(m,n)$, совпадающий с мощностью носителя соответствующей булевой
функции, равен $2^{n-m}$.

\begin{remark}\label{rem:RM3} Отметим, что аналогичным образом множество элементов пространства
$\mathbb{V}_{n,k}({\mathrm{GF}(q)})$  можно рассматривать как линейный код
длины $k^n$, мощности $q^{(k-1)^n}$ и с кодовым расстоянием $2^n$. 
В
частности,  унитрейды в $Q^n_3$ образуют двоичный код длины $3^n$,
мощности $2^{2^n}$ и с кодовым расстоянием $2^n$.

В случае, когда $k=q$, пространство 
$\mathbb{V}_{n,k}(\mathrm{GF}(q))$ имеет достаточно естественное представление в терминах полиномов:
оно состоит из всех функций, которые ортогональны 
любому моному, не зависящему существенно хотя бы от одной из $n$ переменных. Легко понять, что базисом 
$\mathbb{V}_{n,q}({\mathrm{GF}(q)})$ является множество всех мономов, 
у которых степень каждой переменной не превосходит $q-2$. Таким образом, 
$\mathbb{V}_{n,q}(\mathrm{GF}(q))$ можно рассматривать как один из вариантов недвоичного обобщения кодов Рида--Маллера. В частности, один из результатов \S\,\ref{ss:ter-lin}, следствие~\ref{c:3}, можно трактовать в терминах весового распределения этого кода при $q=3$: каждая третья компонента этого распределения нулевая.
\end{remark}

В \cite{Pot12} указана связь между троичными битрейдами и
равномерностью распределения единиц булевой функции по граням.
Булева функция называется {\it почти уравновешенной в гранях}\, если
число нулей и единиц функции отличается не более чем на $2$ в любой
грани любого размера.

\bpro\label{proPot12} Пусть булева функция $f$  уравновешена в
гранях, а $p(x)=x_1\oplus\cdots\oplus x_n$
--- счётчик чётности. Тогда унитрейд
 $U[f\oplus p]$ является битрейдом.  \epro

Из утверждения \ref{probit00} следует, что если булева функция $f$
соответствует битрейду $U[f]$, то и её подфункции, полученные
фиксацией некоторых переменных, также соответствуют битрейдам в
гиперкубах меньшей размерности.

\section{Мощностной спектр множества битрейдов}\label{s:spectr}

В этом разделе будут доказаны свойства спектра мощностей троичных битрейдов, а также битрейдов и унитрейдов малой мощности в произвольных гиперкубах.

\subsection{Мощности унитрейдов и весовой спектр кодов Рида--Маллера}\label{ss:rm}

Из утверждения \ref{probit06} следует, что минимальная мощность
непустого унитрейда размерности $n$ равна $2^n$. В \cite{Pot13} было
доказано

\bpro \label{probit101} Любой унитрейд
 $U\subset Q_k^n$, мощность которого удовлетворяет неравенствам
$2^{n+1}>|U|\geq 2^n$, является битрейдом, имеет $\rank(U)=2$ и
мощность $|U|=2^{n+1}- 2^{s+1} $, где $s\in \{0,\dots,n-1\}$.\epro

Используя результаты исследований весового спектра кодов Рида--Маллера,
можно сильно сузить спектр гипотетических малых (от $2^{n}$ до $5\cdot2^{n-1}$) мощностей
унитрейдов размерности $n$.

\bpro \label{probit01} 
Пусть $U$ --- унитрейд в $Q_k^n$, $k=2^\tau$.
Тогда
существует вектор $u\in \mathcal{R}((\tau-1)n,\tau n)$ такой, что $|U|=\mathrm{wt}(u)$.
\epro 
\proofr Пусть $f=\chi_{_U}:Q_k^n\rightarrow \{0,1\}$.
Рассмотрим произвольное взаимно однозначное отображение $\psi:Q_2^\tau\rightarrow
Q_k$.   Пусть булева функция $F$ определена равенством
$$F=f(\psi(x_1,...,x_\tau),
\psi(x_{\tau+1},x_{2\tau}),
\dots,
\psi(x_{\tau(n-1)+1},x_{\tau n})).$$
Проверим, что $\deg F\leq (\tau-1)n$. 
Рассмотрим любую грань $\Delta$ булева
гиперкуба размерности $(\tau-1)n+1$.
Поскольку  $\Delta$ получена фиксацией значений $n-1$ переменных, найдётся 
$i$ от $0$ до $n-1$ такое, что значения 
переменных 
$x_{\tau i+1}$, \ldots, $x_{\tau i+\tau}$ 
не фиксированы в грани $\Delta$. Тогда по
определению унитрейда  при фиксированных значениях всех остальных
переменных кроме $x_{\tau i+1}$, \ldots, $x_{\tau i+\tau}$ функция $F$ принимает значение
$1$ чётное число раз. Значит она принимает значение $1$ чётное число
раз в грани $\Delta$. Из утверждения \ref{stsp1} следует, что $\deg
F\leq (\tau-1)n$. Тогда набор значений функции $F$ содержится в коде
$\mathcal{R}((\tau-1)n,\tau n)$. Следовательно, унитрейду $U$ соответствует
некоторый вектор $u\in \mathcal{R}((\tau-1)n,\tau n)$. А поскольку функции $F$ и
$f$ одинаковое число раз принимают значение $1$, имеем $|U|=\mathrm{wt}(u)$.
\proofend

\begin{remark}\label{t-trade}
 Аналогично утверждению \ref{probit01} доказывается следующий
 факт: 
 битрейду в $Q_k^n$, $k=2^\tau$ соответствует $[t]$-трейд
 в $Q_2^{\tau n}$, где $t=n-1$ и $[t]$-трейд в $Q_k^{N}$ определяются
 как пара непересекающихся множеств вершин из $Q_k^{N}$,
 разность характеристических функций которых имеет сумму $0$
 по любой $(N-t)$-мерной грани. 
 $[t]$-трейды естественным образом соответствуют разностям ортогональных массивов и алгебраических $t$-дизайнов в графах Хэмминга.
 Аналогичное нашему исследование мощностей малых двоичных
 $[t]$-трейдов было проведено в работе \cite{GKKK}.
 В частности, построены трейды с мощностями
 из серий, рассмотренных нами в \S\,\ref{ss:spectr+}.
\end{remark}

В \cite{KT} и \cite{KTA}  показано, что ненулевые вершины кода
$\mathcal{R}(m,n)$ могут иметь вес только вида $\alpha_m2^{n-m}$,
где $\alpha_m=2-2^{-k}$, $k=0,\dots,n-m-1$, или
$\alpha_m=2+2^{-k}$, $k=2,\dots,\lfloor \frac{n-m}2 \rfloor$, или
$\alpha_m=2+\frac12-2^{-k}$, $k=1,\dots,n-m-1$, или
$\alpha_m=2+\frac12-2^{-k}-2^{-(k+1)}$, $k=3,\dots,n-m-2$, или
$\alpha_m\geq 2+\frac12$.

\btheorem Мощность унитрейда в $Q^n_k$ может принимать значения только вида вида
$\alpha_n2^{n}$, где \\ $\alpha_n=2-2^{-k}$, $k=0,\dots,n-1$, или \\
$\alpha_n=2+2^{-k}$, $k=2,\dots,\lfloor \frac{n}2 \rfloor$, или \\
$\alpha_n=2+\frac12-2^{-k}$, $k=1,\dots,n-1$, или \\
$\alpha_n=2+\frac12-2^{-k}-2^{-(k+1)}$ $k=3,\dots,n-2$, или \\
$\alpha_n\geq 2+\frac12$.  \etheorem 
\proofr Любой унитрейд в $Q^n_k$
имеет мощность не меньше $2^n$ (утверждение \ref{probit06}). Поэтому
унитрейд, пересекающийся с $5$ гипергранями некоторого направления,
имеет мощность не менее $5\cdot2^{n-1}$ (последний случай в утверждении теоремы). Если же унитрейд
пересекается не более чем с $4$ гипергранями любого направления и
является подмножеством вершин подграфа, изоморфного $Q^n_4$. 
В этом случае 
требуемое следует из утверждения~\ref{probit01} ($\tau=2$) и перечисленных выше
результатов из~\cite{KT} и~\cite{KTA}. \proofend

\subsection{Троичные битрейды как линейные функции над полем GF(3)}\label{ss:ter-lin}

Вначале выберем подходящий базис в пространстве
$\mathbb{V}_{n,3}({\mathrm{GF}(3)})$. Здесь нам будет удобно считать, что
$Q_3=\{0,\pm 1\}=\mathrm{GF}(3)$. Пусть $s_0(a)=1$ и $s_1(a)=a$. Определим
функции $s_\alpha: Q^n_3 \rightarrow Q_3$, $\alpha\in Q^n_2$,
равенствами $s_\alpha(x)=s_{\alpha_1}(x_1)\cdots s_{\alpha_n}(x_n)$.

\bpro\label{probit07} (a) $s_\alpha \in \mathbb{V}_{n,3}({\mathrm{GF}(3)})$;\\
(b) $\{s_\alpha \ |\ \alpha\in Q^n_2\}$ --- базис в
$\mathbb{V}_{n,3}({\mathrm{GF}(3)})$;\\
с) $\langle s_\alpha, s_\beta\rangle_3 = \sum\limits_{x\in Q^n_3}
s_\alpha(x) s_\beta(x)=0$, за исключением случая
$\alpha=\beta=\overline{1}$. \epro

\proofr
(a) Следует из того, что $\sum\limits_{a\in Q_3}
s_0(a)=\sum\limits_{a\in Q_3} s_1(a)=0$.

(b) Заметим, что $\alpha\in \supp(s_\alpha)$, а  $\alpha \in
\supp(s_\beta)$, только если $\beta\preceq \alpha$,
$\alpha,\beta\in Q^n_2$. Упорядочим функции $s_\alpha$ в порядке
(частичном) убывания от $\overline{1}$ до $\overline{0}$. Носитель
каждой следующей функции содержит точку, которая не содержится в
носителях предыдущих функций. Поэтому на каждом шаге мы будем
получать линейно независимое семейство функций. Как показано выше,
$\dim (\mathbb{V}_{n,3}({\mathrm{GF}(3)}))=2^n$, следовательно, семейство
линейно независимых функций мощности $2^n$ является базисом.

(c) Пусть набор $\alpha$ имеет нулевую координату. Без ограничения
общности будем считать,
что $\alpha_n=0$. Тогда справедливы равенства\\
$\sum\limits_{x\in Q^n_3} s_\alpha(x) s_\beta(x)=\sum
s_{\alpha_1}(x_1)s_{\beta_1}(x_1)\cdots s_{\alpha_{n-1}}(x_{n-1})
s_{\beta_{n-1}}(x_{n-1})\sum\limits_{x_n\in Q_3} s_0(x_n)
s_{\beta_{n}}(x_{n})=$\\
$\sum s_{\alpha_1}(x_1)s_{\beta_1}(x_1)\cdots
s_{\alpha_{n-1}}(x_{n-1})
s_{\beta_{n-1}}(x_{n-1})\sum\limits_{x_n\in Q_3}
s_{\beta_{n}}(x_{n})=0$. $\blacktriangle$

\bcorol\label{c:3}
Для любого $f\in\mathbb{V}_{n,3}({\mathrm{GF}(3)})$ верно $\langle f,
f \rangle_3\equiv 0, 2^n\bmod 3$.
\ecorol

\proofr Из утверждения \ref{probit07} следует, что\\
$f=\sum\limits_{\alpha\in Q^n_2} a_\alpha s_\alpha$ и  $\langle f, f
\rangle_3= \sum\limits_{\alpha,\beta\in Q^n_2} a_\alpha a_\beta
\sum\limits_{x\in Q_3^n} s_\alpha(x) s_\beta(x)= a_{\overline{1}}^2
\sum\limits_{x\in Q_3^n}s_{\overline{1}}^2=a_{\overline{1}}^2
|\supp(s_{\overline{1}})|= a_{\overline{1}}^22^n$, где все операции
выполняются в поле $\mathrm{GF}(3)$. $\blacktriangle$

Как отмечено в замечании~\ref{rem:RM3}, пространство 
$\mathbb{V}_{n,3}({\mathrm{GF}(3)})$ является троичным кодом типа 
Рида--Маллера, а следствие~\ref{c:3} означает, что
весовой спектр этого кода имеет нули в каждой третьей компоненте.

\btheorem \label{probit102} Для любого битрейда $B\subset Q_3^n$
верно, что $|B|\equiv 0, 2^n\bmod 3$. 
\etheorem

\proofr Рассмотрим функцию $b:Q_3^n\rightarrow Q_3$, принимающую
значения $1$ и $-1$ на двух долях битрейда $B$ и $0$ в остальных
точках. Тогда $b\in\mathbb{V}_{n,3}({\mathrm{GF}(3)})$ и $|B|=\langle b, b
\rangle_3\equiv 0, 2^n\bmod 3$. \proofend

\bcorol\label{probit1020} 
Пусть $U\subset Q_3^n$ --- унитрейд мощности $2^{n+1}$. Тогда $U$ не битрейд. \ecorol

\subsection{Некоторые свойства мощностного спектра троичных
битрейдов}\label{ss:some}

Рассмотрим несколько простых свойств битрейдов. Пусть $f$ --- булева
функция от переменных $x_1,\dots,x_n$, а $g$ --- булева функция от
переменных $y_1,\dots,y_m$ и наборы переменных не пересекаются.
Обозначим, через $f(x)g(y)$ булеву функцию от $n+m$ переменных. В
\cite{Pot13} имеется

\bpro\label{probit20} Если $U[f]\subset  Q_3^n$ и $U[g]\subset
Q_3^m$
--- битрейды, то $U[f(x)g(y)]\subset Q_3^{n+m}$ --- битрейд и
$|U[f(x)g(y)]|=|U[f]||U[g]|$.\epro

В частности, в качестве $g$ можно рассмотреть линейную функцию
$g(y)=\bigoplus y_i$ и  как следствие получить свойство: если в
$Q_3^n$ имеется битрейд мощности $a$, то в $Q_3^{n+m}$ имеется
битрейд мощности $a2^m$. Эту конструкцию можно рассматривать как
частный случай декартова произведения двух битрейдов. Декартово
произведение двух  битрейдов является битрейдом не только в троичном
гиперкубе, но и в гиперкубах над произвольным алфавитом. Унитрейды,
которые можно представить в виде декартова произведения, будем
называть {\it разложимыми}.

\btheorem[о построении редуцируемых битрейдов]
\label{probit1021} Пусть $B\subset Q_k^n$, $C\subset Q_k^n$  --- битрейды. 
Тогда в $Q_k^{n+m}$ имеются битрейды мощности 
$|B|\cdot|C|$, $2^m|B|$,  $k^m|B|$.
\etheorem
\proofr
Битрейды мощности $|B|\cdot|C|$ и $2^m|B|$ строются декартовым произведением,
последнее утверждение следует по индукции из случая $m=1$, который разберем отдельно.
Пусть функция 
$b:Q_k^n\rightarrow \{-1,0,1\}$ принимает значения $1$ и
$-1$ на двух долях некоторого битрейда и $0$ в остальных точках.
Тогда функция $b':Q_k^{n+1}\rightarrow \{-1,0,1\}$, заданная равенством
$b'(x_1,\dots,x_n,x_{n+1})=b(x_1,\dots,x_{n-1},(x_n+x_{n+1})\bmod
k)$, определяет битрейд размерности $n+1$ и мощности $k|B|$. 
\proofend

\bpro\label{probit21} Если унитрейд $U\subset  Q_3^n$ имеет пустой
ретракт по некоторому направлению, то он эквивалентен унитрейду
$U[f]$, где $f$
--- булева функция не зависящая от одной из переменных. Причём $U$
является битрейдом тогда и только тогда, когда битрейдом является
любой из его непустых ретрактов по тому же направлению. \epro
\proofr Без ограничения общности можно считать, что $U\cap
\{x_n=-1\}=\varnothing$. Тогда из определения унитрейда имеем $U\cap
\{x_n=0\}=U\cap \{x_n=1\}=U'$ и $U=U'\times \{0,1\}$. Пусть
$U'=U[f]$ и $g(x_1,\dots,x_{n})=f(x_1,\dots,x_{n-1})$. Тогда
$U=U[g]$. Из утверждений \ref{probit00} и \ref{probit20} следует,
что унитрейды $U$ и $U'$ являются битрейдами одновременно. \proofend

Поскольку каждая одномерная грань гиперкуба пересекается с
унитрейдом не более чем по двум вершинам, троичный унитрейд $U$
размерности $n$ содержит не более $2\cdot3^{n-1}$ вершин. Причём в
случае равенства $|U|=2\cdot3^{n-1}$, дополнение $Q^n_3\setminus U$
является МДР-кодом. Как известно (см., например, \cite{KP09}), в
$Q^n_3$ имеется единственный с точностью до эквивалентности МДР-код,
который является линейным (т.\,е. аффинным подпространством над полем
$\mathrm{GF}(3)$). Поэтому унитрейд максимальной мощности единственный,
что отражено в первой части следующей теоремы.

 \btheorem[о битрейдах большой мощности] (a) В $Q_3^n$ имеется только один с точностью до
эквивалентности унитрейд $B$ максимальной мощности $2\cdot3^{n-1}$,
который является
битрейдом и может быть задан равенством 
$$B= \{x\in Q_3^n \ |\
x_1+\dots+x_n\not \equiv 0 \, \bmod\, 3\}=U[\ell],
\qquad \ell(x_1,\dots,x_n)= \bigoplus_{i<j}x_ix_j.$$ 

(b) В $Q_3^{n}$ имеются битрейды мощности $14\cdot3^{n-3}$.

(c) В $Q_3^n$ не существует битрейдов мощности, 
промежуточной
между $14\cdot3^{n-3}$ и $2\cdot3^{n-1}$.  
\etheorem
 \proofr
Осталось доказать (b) и (c).
В $Q_3^3$ имеется битрейд $U[x_1\oplus x_2\oplus x_3
\oplus 1]$ мощности $14$. 
По теореме~\ref{probit1021} имеем (b).

Докажем (c) по индукции.
Базой индукции является случай $n=3$, который легко проверить непосредственно.
 По индукции. Пусть $U^{n+1}\subset Q_3^{n+1}$ --- некоторый битрейд.
Если три различных ретракта  какого-то одного направления имеют
мощности меньше $2\cdot3^{n-1}$, то $|U^{n+1}|\leq 14\cdot3^{n-2}$
по предположению индукции. Предположим  в $U^{n+1}$ найдутся по
одному ретракту каждого направления мощности $2\cdot3^{n-1}$. Без
ограничения общности можно полагать, что каждый из этих ретрактов
соответствует  нулевому значению некоторой переменной ($x_i=0$). По
предположению индукции они с точностью до эквивалентности заданы
линейными над $\mathrm{GF}(3)$ функциями. Тогда $U^{n+1}=U[f]$, где $f=\ell$
или $f=\ell\oplus x_1\cdots x_{n+1}$. Однако, вторая функция имеет
недвудольный $4$-x мерный ретракт. Поэтому в случае когда
$U[f]\subset Q_3^{n+1}$ --- битрейд, имеем $f=\ell$ и
$|U^{n+1}|=|U[f]|=2\cdot 3^n$.
 \proofend

\subsection{Битрейды и расстояния между мономами}\label{ss:monom}

Пусть $B\in Q^n_k$ --- битрейд. В \S\,\ref{s:lin} была определена функция
$b[B]:Q^n_k\rightarrow \{0,\pm 1\}$, которая принимает значение $1$
на одной доле битрейда $B$, $-1$ на другой доле битрейда $B$ и $0$ в
остальных вершинах. В этом разделе мы будем рассматривать $b[B]$ как
функцию действующую в $\mathbb{R}$, т.\,е. $b[B]\in
\mathbb{V}_{n,3}(\mathbb{R})$.  В случае $k=3$ такую функцию
можно определить ровно двумя способами $b[B]$ и $-b[B]$, поскольку
граф $\Gamma B$ связен. В дальнейшем мы рассматриваем только
такой случай и обычно не уточняем, каким из двух способов выбран
знак функции $b[B]$. Для краткости введём обозначения
$b_v=b[U[x^v]]$ и $b_V=b[U[f^V]]$, где  $v\in Q^n_3$, $V\subset
Q^n_3$.

\bpro Пусть $B, B'\in Q^n_3$ --- битрейды и $b[B]b[B']\neq 1$. Тогда
$S=\supp (b[B]+b[B'])$ --- битрейд и $b[S]=b[B]+b[B']$.\epro \proofr
Если $b[B], b[B']\in \mathbb{V}_{n,3}(\mathbb{R})$, то
$b[B]+b[B']\in \mathbb{V}_{n,3}(\mathbb{R})$. По условию
$b[B]b[B']\neq 1$, следовательно, $(b[B]+b[B'])(Q^n_3) \subseteq
\{0,\pm 1\}$. Тогда по определению $S=\supp (b[B]+b[B'])$ ---
битрейд. \proofend

Будем говорить, что функции $b_v$ и $b_u$ {\it согласованы}, если
$b_ub_v\neq 1$.

\bpro\label{probit22} Пусть $v, u\in Q^n_3$.  Если найдётся вершина
$x\in Q^n_3$ такая, что $b_v(x)b_u(x)=-1$, то пара функций $b_u$ и
$b_v$ согласована. \epro \proofr
 Пересечение битрейдов $U[x^u]$ и $U[x^v]$ является булевым
подкубом в  грани, соответствующей совпадающим координатам наборов
$u$ и $v$. Из связности пересечения $U[x^u]\cap U[x^v]$ следует, что
унитрейд $U[x^u\oplus x^v]$ является битрейдом. Как было замечено
выше, любой битрейд разделяется на доли единственным образом.
\proofend

\bcorol\label{probit222} Любой унитрейд ранга $2$ является битрейдом.
\ecorol

\bpro\label{probit25} Унитрейд $U[f^V]$ является битрейдом, если для
каждого $v\in V$ можно выбрать функции (знак функции) $b_v$ так,
чтобы функция $g=\sum\limits_{v\in V}b_v$ принимала значения только
из множества $\{0,\pm 1\}$. \epro \proofr Поскольку $b_v\in
\mathbb{V}_{n,3}(\mathbb{R})$ для любого $v\in V$, то и $g\in
\mathbb{V}_{n,3}(\mathbb{R})$. Тогда из условия следует, что $\supp
(g)$ есть унитрейд. Очевидно $U[f^V]\subset \supp (g) $. Тогда
$\supp (g) = U[f^V]$, т.\,е. $g=\pm b_V$. \proofend

\bpro\label{bitrade90} Пусть $V\subset Q^n_3$ и  $|V|=3$ (т.\,е.
$\rank(U[f^V])\leq 3$). Если все вершины множества $V$ различаются
только в двух координатах, либо две вершины множества $V$
различаются только в одной координате, то $U[f^V]$
--- битрейд. \epro \proofr
 Пусть
$V=\{u,v,w\}$ и вершины множества $V$ различаются только в двух
координатах. Тогда с помощью утверждения \ref{probit21} можно
перейти к двумерному случаю, когда все унитрейды являются
битрейдами.

 Если $d(v,u)=1$, то $x^v\oplus x^u=x^o$, т.\,е.
$\rank(U[f^V])= 2$ и требуемое следует из следствия \ref{probit222}.
\proofend

Будем говорить, что три вершины $\{u,v,w\}\subset Q^n_3$ {\it
находятся в общем положении}, если найдётся координата, в которой
они попарно различаются. В этом случае найдутся точки $a,a',a''\in
Q_3^{n}$,  в которых носители битрейдов $b_u,b_v,b_w$ пересекаются
только попарно, т.\,е.  $a\in (\supp(b_{v})\cap\supp(b_{u}))\setminus
\supp(b_{w})$, $a'\in (\supp(b_{v})\cap\supp(b_{w}))\setminus
\supp(b_{u})$, $a''\in (\supp(b_{w})\cap\supp(b_{u}))\setminus
\supp(b_{v})$.

Вначале рассмотрим унитрейды $U[f^V]$ ранга $3$, когда три вершины
$V=\{u,v,w\}$ не находятся в общем положении.

\bpro\label{bitrade91} Пусть $V=\{u,v,w\}\subset Q^n_3$.\\
(a) Если любая координата набора $w$ совпадает с соответствующей
координатой набора $u$ или набора $v$, то $U[f^V]$
--- битрейд.\\
(b) Если любая координата на наборах $u$, $v$, $w$ принимает не более двух
значений и не выполнено условие (a), то $U[f^V]$ не является
битрейдом. \epro 

\proofr
(a) Рассмотрим случай, когда нет координаты, в которой все три
вершины $u$, $v$, $w$ совпадают. Без ограничения общности можно считать,
что $u=\overline{0}$, $v=\overline{1}$, $w\in \{0,1\}^n$. По
утверждению \ref{probit22} функции $b_v$ и $b_u$ можно выбрать так
чтобы вещественная сумма $b_v+b_w$ и $b_u+b_w$ принимала значения
только из множества $\{0,\pm 1\}$. Из условия видно, что
$U[x^{\overline{0}}]\cap U[x^{\overline{1}}] =\{
\overline{1}\}\subset U[x^w]$. Поэтому $(b_v+b_w)(\overline{0})=0$.
Следовательно, функция $b_v+ b_u+b_w$ принимает значения только  из
множества $\{0,\pm 1\}$. Требуемое следует из утверждения
\ref{probit25}.

Случай, когда все три вершины $u$, $v$, $w$ совпадают в некоторой
координате сводится к рассмотренному с помощью утверждения
\ref{probit21}.

(b) Если любая тройка координат в наборах $u$, $v$, $w$ удовлетворяет
условию (a), то наборы $u$, $v$, $w$ также удовлетворяют этому условию.
 Без ограничения общности можно считать, что условию (a) не удовлетворяют первые тройки
координат наборов $u$, $v$, $w$ и они равны соответственно $e_1=(1,0,0)$,
$e_2=(0,1,0)$ и $e_3=(0,0,1)$. Поскольку в каждой координате наборы
$u$, $v$, $w$ принимают только два значения, то найдётся трёхмерная грань,
на которой $x^v\oplus x^u\oplus x^w= x_1\oplus x_2\oplus x_3= f'$.
Непосредственная проверка показывает, что $U[f']$ не является
битрейдом. Тогда из утверждения~\ref{probit00} следует, что и
$U[f^V]$ не является битрейдом. \proofend

Теперь рассмотрим унитрейды $U[f^V]$ ранга $3$, когда три вершины
$V=\{u,v,w\}$  находятся в общем положении.

\bpro\label{bitrade9} Пусть $V=\{u,v,w\}\subset Q^n_3$.\\
(a) Если сумма попарных расстояний между вершинами множества $V$
нечётная, то они находятся в общем положении.\\
(b) Если вершины множества $V$ находятся в общем положении, то  из
согласованности пары $b_v$, $b_u$ и пары $b_v$, $b_w$ следует
согласованность пары $b_u$, $b_w$ тогда и только тогда, когда  сумма
попарных расстояний между вершинами множества $V$ нечётная.  \epro
\proofr Если вершины $u$, $v$, $w$ не находятся в общем положении, то
каждая координата даёт вклад $0$ или $2$ в сумму
$d(u,v)+d(v,w)+d(u,v)$. Поэтому в этом случае сумма расстояний
чётная. Пункт (a) доказан.

 Рассмотрим случай, когда
$V=\{\overline{0},\overline{1},\overline{-1}\}$. Нетрудно видеть,
что унитрейды $U[x^{\overline{0}}]$ и $U[x^{\overline{1}}]$
пересекаются по одной точке $\overline{1}$. Разделим битрейд
$U[x^{\overline{0}}\oplus x^{\overline{1}}]$ на две доли. Вершины
той же чётности, что и $\overline{1}$ в гиперкубах
$U[x^{\overline{0}}]$ и $U[x^{\overline{1}}]$ должны принадлежать
разным долям. Следовательно, вершины $\overline{0}$ и
$\overline{-1}$ оказываются в разных долях. Рассмотрим гиперкуб
$U[x^{\overline{-1}}]$. Ясно, что вершины $\overline{0}$ и
$\overline{-1}$ оказываются в разных долях если и только если $n$
нечётно. Тогда из согласованности пары $b_{\overline{0}}$,
$b_{\overline{1}}$ и пары $b_{\overline{0}}$, $b_{\overline{-1}}$
 следует согласованность пары $b_{\overline{1}}$,
$b_{\overline{-1}}$ при нечётном $n$ и несогласованность пары
$b_{\overline{1}}$, $b_{\overline{-1}}$ при чётном $n$.

Проведём дальнейшее доказательство по индукции. Предположим, что при
$n-1$ утверждение доказано.

 Пусть вершины $u$, $v$, $w$ попарно различаются  в каждой координате.
Тогда утверждение следует из рассмотренного выше случая. В противном
случае найдётся координата, в которой не все вершины $u$, $v$, $w$ попарно
различаются. После удаления этой координаты укороченные наборы
$v'$, $u'$, $w'$ также находятся в общем положении. Без ограничения
общности можно полагать, что удалённая координата была последней.
Очевидно, $d(u',v')+d(v',w')+d(u',v')=d(u,v)+d(v,w)+d(u,v)$, если
$u_n=v_n=w_n$, и $d(u',v')+d(v',w')+d(u',v')=d(u,v)+d(v,w)+d(u,v)-2$
в противном случае. Тогда утверждение верно для функций $b_{v'}$,
$b_{u'}$ и $b_{w'}$ по предположению индукции. Поскольку последняя
координата не принимает одно из значений $\{0,\pm 1\}$ на наборах
$u$, $v$, $w$, то найдётся гипергрань $x_n=\delta$ по последнему
направлению такая, что $b_{v}(x\delta)=b_{v'}(x)$,
$b_{u}(x\delta)=b_{u'}(x)$, $b_{w}(x\delta)=b_{w'}(x)$  для всех
$x\in Q^{n-1}_3$. Тогда  для функций $b_{v}$, $b_{u}$ и $b_{w}$
требуемое следует из утверждения \ref{probit22}. \proofend

\bcorol\label{bitrade97}  Если сумма попарных расстояний между
вершинами множества $V=\{u,v,w\}\subset Q^n_3$ нечётная, то $U[f^V]$
--- битрейд.
\ecorol

\bcorol\label{bitrade98}
 Если вершины множества $V=\{u,v,w\}\subset Q^n_3$ находятся в общем
 положении, различаются не менее чем в трёх координатах,
  никакие две из них не находятся на расстоянии $1$ и
 сумма попарных расстояний между вершинами
множества $V$ чётная, то $U[f^V]$ не является битрейдом. \ecorol
\proofr из предложения \ref{bitrade9} следует, что набор функций
$b_u$, $b_v$, $b_w$ не может быть попарно согласованным. Из условия:
вершины $u$, $v$, $w$  различаются не менее чем в трёх координатах и
никакие две из них не находятся на расстоянии $1$ следует, что
множества $\supp(b_{w})\setminus((\supp(b_{v})\cup\supp(b_{u})) )$,
$\supp(b_{u})\setminus((\supp(b_{v})\cup\supp(b_{w})) )$,
$\supp(b_{v})\setminus((\supp(b_{w})\cup\supp(b_{u})) )$ являются
связными. Тогда из несогласованности следует недвудольность
унитрейда. \proofend

\bpro Если все попарные расстояния между различными точками
множества $V\subset Q^n_3$ нечётные, то $U[f^V]$ --- битрейд. \epro

\proofr Покажем по индукции, что можно выбрать функции $b_v$, $v\in
V$, так что все функции $b_v$ попарно согласованы. При $|V|=3$
требуемое следует из утверждения \ref{bitrade9}. Пусть для множеств
$V\subset Q^n_3$, $|V|=k$, утверждение верно.  Рассмотрим $u\not \in
V$, находящуюся на нечётном расстоянии от любой точки из $V$.
Выберем функцию $b_u$ согласовано с $b_v$ для некоторого $v\in V$.
Тогда из утверждения \ref{bitrade9} следует, что $b_u$ согласовано с
$b_w$ для любого $w\in V$. Таким образом, утверждение доказано при
$|V|=k+1$.

Если все функции $b_v$, $v\in V$, согласованы, то $\sum\limits_{v\in
V}b_v$ принимает значения только из множества $\{0,\pm 1\}$. Тогда
$U[f^V]$ --- битрейд по утверждению \ref{probit25}. \proofend

\subsection{Вычислительные результаты}\label{ss:tables}
В этом разделе мы приведем результаты вычислений числа битрейдов при помощи ЭВМ.
Удалось вычислить число различных битрейдов до размерности $n=7$, а число неэквивалентных --- до размерности $n=6$. 
Метод перечисления не сильно отличается от прямого перебора, 
поэтому мы не будем описывать алгоритм в подробностях.
Для перечисления всех битрейдов в $Q_3^n$,
в качестве одного из ретрактов подставлялись по одному представителю из каждого 
из $N'(n-1)$ классов эквивалентности функций, найденных на предыдущем шаге.
После этого для параллельного ретракта проводился 
поточечный перебор значений функции,
с очевидной проверкой выполнимости условия на сумму по одномерной грани.
Полуенное число решений умножалось на число представителей в классе эквивалентности первого ретракта.
Вычисление $N(7)$ заняло два года процессорного времени (в рассчете на одно ядро процессора). 
Результаты вычислений до $n=6$ проверены при помощи следующей техники
двойного подсчета (см. \cite{KO:alg}):
мощность каждого класса эквивалентности, вычисленная через мощность группы автоморфизмов его представителя, совпадает с числом представителей, 
найденных в процессе полного перебора.
В следующей таблице $N(n)$ --- число различных (включая тождественно нулевую) 
$\{-1,0,1\}$-функций на $Q_3^n$, 
у которых сумма по каждой одномерной грани равна $0$ 
(то есть все три значения в грани либо нулевые, либо попарно различные); 
$N'(n)$ --- число неэквивалентных таких функций. Последняя колонка отражает среднюю половинную мощность (число $-1$-ц) битрейда (в скобках приведено среднеквадратичное отклонение); из этой статистики исключен пустой трейд, так как он не удовлетворяет критерию хи-квадрат.

$$
\begin{array}{r|l|l|l|l|l|l}
n& N'(n) & N(n) & \ln N(n) & \ln \ln N(n)&    \\ \hline
0& 2           & 3                   &  1.098 & 0.094       &  \\
1& 2           & 7                   &  1.945 & 0.665  (+0.571) & 1  \\
2& 3           & 31                  &  3.433 & 1.233  (+0.567) & 2.4 \ (\pm0.490) \\
3& 5           & 403                 &  5.998 & 1.791  (+0.557) &  6.448 =2.539^2\ (\pm 1.188)  \\
4& 13    & 29875    & 10.304 & 2.332  (+0.541) & 17.960=2.619^3\ (\pm2.342) \\
5& 92    & 32184151 & 17.286 & 2.849  (+0.517) & 50.527=2.666^4\ (\pm4.776) \\
6& 25493       & 1488159817231       & 28.028 & 3.333 (+0.483) &  142.25=2.695^5 \ (\pm10.07) \\
7& >2187260868 & 6171914027409468739 & 43.266 & 3.767 (+0.434)
& 398.17=2.712^6\ (\pm22.59)
\end{array}
$$
Далее мы выпишем распределение битрейдов по мощностям. 
Для каждого $n$ указано число битрейдов (именно двудольных унитрейдов, то есть число функций будет вдвое больше) мощности $2^n$, $2^n+2$, $2^n+4$, \ldots, $2\cdot 3^{n-1}$.

$n=1$: 3.

$n=2$: 9, 6.

$n=3$: 27, 0, 54, 108, 0, 12.

$n=4$: 81, 0, 0, 0, 324, 0, 1296, 648, 0, 3888, 2844, 0, 4536, 1296, 0, 0, 0, 0, 0, 24.

$n=5$: 243, 0, 0, 0, 0, 0, 0, 0, 
1620, 0, 0, 0, 9720, 0, 9720, 3888, 0, 0, 58320, 0, 41580, 77760, 0, 116640, 
301320, 0, 259200, 660960, 0, 480816, 1368576, 0, 1156680, 2468880, 0, 1415232, 
2721600, 0, 1148040, 2185056, 0, 583200, 816480, 0, 90720, 104976, 0, 10800, 0, 0,
0, 0, 0, 0, 0, 0, 0, 0, 0, 0, 0, 0, 0, 0, 0, 48.

$n=6$: 729, 0, 0, 0, 0, 0, 0, 0, 0, 0, 0, 0, 0, 0, 0, 0, 7290, 0, 0,
0, 0, 0, 0, 0, 58320, 0, 0, 0, 87480, 0, 69984, 23328, 0, 0, 0, 0, 524880, 0, 0,
0, 370980, 0, 1399680, 0, 0, 699840, 2099520, 0, 4811400, 466560, 0, 2799360,
7290000, 0, 16562880, 2099520, 0, 6998400, 15244848, 0, 49968576, 19012320, 0,
46889280, 48114000, 0, 149999040, 48988800, 0, 173560320, 158793696, 0, 431451360,
203303520, 0, 593464320, 402077520, 0, 1226726208, 655983360, 0, 1759957632,
1275108480, 0, 3455693280, 1610681760, 0, 4922674560, 3332579760, 0, 8667868320,
4840793280, 0, 12263996160, 7124630400, 0, 19261521360, 10458292320, 0,
25982259840, 15546805632, 0, 37437240960, 17859890880, 0, 44159904000,
26492909760, 0, 56014493760, 28054486080, 0, 58200653952, 29447634240, 0,
63563901120, 30536701920, 0, 53914973760, 27520508160, 0, 44905723488,
18151205760, 0, 28971976320, 13573863360, 0, 17778852000, 6267067200, 0,
7903992960, 2932269264, 0, 2917632960, 1190894400, 0, 772623360, 243544320, 0,
299531520, 100077120, 0, 33592320, 41290560, 0, 5598720, 6298560, 0, 0, 1179360,
0, 3079296, 3429216, 0, 0, 0, 0, 0, 77760, 0, 0, 0, 0, 0, 0, 0, 0, 0, 0, 0, 0, 0,
0, 0, 0, 0, 0, 0, 0, 0, 0, 0, 0, 0, 0, 0, 0, 0, 0, 0, 0, 0, 0, 0, 0, 0, 0, 0, 0,
0, 0, 0, 0, 0, 0, 0, 0, 0, 0, 0, 0, 0, 96.

$n=7$: { \footnotesize 2187, 0, 0, 0, 0, 0, 0, 0, 0, 0, 0, 0, 0, 0,
0, 0, 0, 0, 0, 0, 0, 0, 0, 0, 0, 0, 0, 0, 0, 0, 0, 0, 30618, 0, 0, 0, 0, 0, 0, 0, 
0, 0, 0, 0, 0, 0, 0, 0, 306180, 0, 0, 0, 0, 0, 0, 0, 612360, 0, 0, 0, 734832, 0, 489888, 139968, 0, 0, 0, 0, 0, 0, 0, 0, 3674160, 0, 0, 0, 0, 0, 0, 0, 2585520, 0,
0, 0, 14696640, 0, 0, 0, 0, 0, 14696640, 0, 22044960, 5878656, 0, 0, 48376440, 0, 9797760, 0, 0, 0, 58786560, 0, 105325920, 9797760,  0, 29393280, 252292320, 0,
48988800, 19595520, 0, 0, 264539520, 0, 258660864, 39191040, 0, 58786560,
849465792, 0, 417384576, 64665216, 0, 118599552, 1102248000, 0, 1026927720,
440899200, 0, 117573120, 3241833840, 0, 1646023680, 881798400, 0, 411505920,
5472048960, 0, 4595639328, 1293304320, 0, 1175731200, 13190234400, 0, 6642881280,
4967464320, 0, 2792361600, 23992754688, 0, 14767655616, 8897345856, 0, 5232003840,
48779617824, 0, 29422673280, 19428958080, 0, 13418032320, 86712135552, 0,
56942131680, 43070952960, 0, 29628426240, 174586285440, 0, 108140326560,
100338860160, 0, 60767667072, 333349188480, 0, 211194390960, 197052549120, 0,
133406300160, 633974103504, 0, 394756649280, 437915782080, 0, 284879669760,
1184769633600, 0, 732195313920, 918265662720, 0, 553710608640, 2237431905072, 0,
1396486980000, 1839642114240, 0, 1166266563840, 4133841505920, 0, 2463586985664,
3699825232320, 0, 2436508916352, 7740828063360, 0, 4633391621376, 7373480647680,
0, 4422630481920, 14095095607296, 0, 9008038009152, 14256009258624, 0,
8787284896320, 25936010563680, 0, 15056472533760, 27127920314880, 0,
17656455116160, 46489296385728, 0, 28595219380560, 51865158776256, 0,
30807273127680, 84076043325120, 0, 53419783072320, 96717020198400, 0,
59173324616448, 151123241747520, 0, 87927292938240, 177969306401280, 0,
112318130615040, 267008385528864, 0, 160800706809792, 321325302945024, 0,
188053645282560, 466481709832320, 0, 291050421480000, 570238841894400, 0,
343218104712000, 806560729743456, 0, 464604287555520, 996450456984192, 0,
620225396784000, 1373638417522560, 0, 816009087669552, 1708948701149760, 0,
982341376894080, 2305937678879520, 0, 1414154728014336, 2868447898270080, 0,
1691741075976000, 3783164938709184, 0, 2145564791750016, 4708186964144640, 0, 
2866434717250944, 6089563693805760, 0, 3543718526510400, 7547844801104640, 0,
4235506218558720, 9574787807964288, 0, 5748041164944960, 11803431065489280, 0, 
6789801853310400, 14687203311950400, 0, 8100416661309504, 17888868241614720, 0,
10588066580077056, 21840361606638720, 0, 12363491876450400, 26269417255213440, 0,
14285288100787200, 31462082237108160, 0, 18299032850230272, 37233543050766720, 0,
20709927455180544, 43720111582963200, 0, 23313370728464640, 50763763455713280, 0,
28924064989464960, 58412989843236000, 0, 31813638206316480, 66303578484210816, 0,
34565462372004480, 74583100623265920, 0, 41520901293714528, 82545422286942720, 0,
43826073580183104, 90590698165121280, 0, 45970461276122880, 97356175759100160, 0,
52673395505996160, 103853217900781440, 0, 53561231066238336, 108100043912367360,
0, 53128915099827840, 111570541229580480, 0, 58350442395913920,
112097875912081920, 0, 55765902028032000, 111406452367500864, 0,
52657506351233280, 107575556428968768, 0, 53994455165468160, 102377410528901760,
0, 48778189150406400, 94582347577421760, 0, 42669649619928576, 85745409629443200,
0, 40987358888555520, 75337864538158080, 0, 34117558361470080, 64695169565885376, 
0, 27623026142341056, 53722349499008448, 0,  24230050550665344, 43420216983171840,
0, 18493690067425440, 33854862632935680, 0, 13546143606925440, 25584444523483776,
0, 10783563464378880, 18568192776336000, 0, 7350295708661760, 13016846960075520, 
0, 4835426179046400, 8727631849641600, 0, 3393591789458304, 5646375105114240, 0,
2053100209063680, 3458776067268480, 0, 1173620938855680, 2041794442509696, 0,
728670130929216, 1138311067888512, 0, 375048142382592, 609672959421120, 0,
191348768439360, 307095833018880, 0, 99253170374400, 152237238241536, 0,
46874958318720, 69750838786176, 0, 19948807630080, 31988971163520, 0,
10729375110720, 13802819748480, 0, 3921475057920, 6002269439040, 0, 1657927958400,
2388145213440, 0, 1007425278720, 1073677731840, 0, 520179426144, 484107321600, 0,
144614937600, 245492674560, 0, 205312060800, 136090886400, 0, 89050752000,
73071694080, 0, 20222576640, 49468890240, 0, 45500797440, 28217548800, 0,
15872371200, 12433357440, 0, 3174474240, 2821754880, 0, 12580323840, 4938071040,
0, 3853785600, 1058158080, 0, 0, 529079040, 0, 1410877440, 0, 0, 1162667520, 0, 0,
0, 0, 0, 235146240, 0, 0, 264539520, 0, 0, 0, 0, 0, 0, 0, 0, 12700800, 0, 0, 0, 0,
0, 129330432, 120932352, 0, 71197056, 0, 0, 0, 0, 0, 0, 0, 0, 0, 0, 0, 0, 0, 0, 0,
0, 0, 520128, 0, 0, 0, 0, 0, 0, 0, 0, 0, 0, 0, 0, 0, 0, 0, 0, 0, 0, 0, 0, 0, 0, 0,
0, 0, 0, 0, 0, 0, 0, 0, 0, 0, 0, 0, 0, 0, 0, 0, 0, 0, 0, 0, 0, 0, 0, 0, 0, 0, 0,
0, 0, 0, 0, 0, 0, 0, 0, 0, 0, 0, 0, 0, 0, 0, 0, 0, 0, 0, 0, 0, 0, 0, 0, 0, 0, 0,
0, 0, 0, 0, 0, 0, 0, 0, 0, 0, 0, 0, 0, 0, 0, 0, 0, 0, 0, 0, 0, 0, 0, 0, 0, 0, 0,
0, 0, 0, 0, 0, 0, 0, 0, 0, 0, 0, 0, 0, 0, 0, 0, 0, 0, 0, 0, 0, 0, 0, 0, 0, 0, 0,
0, 0, 0, 0, 0, 0, 0, 0, 0, 0, 0, 0, 0, 0, 0, 0, 0, 0, 0, 0, 0, 0, 0, 0, 0, 0, 0,
0, 0, 0, 192.
}

\subsection{Уточнение мощностного спектра битрейдов}\label{ss:spectr+}

\bpro Пусть  $\rank(U)\geq 4$ и $U\subset Q^n_3$ --- неразложимый
унитрейд размерности $n$. Тогда $|U|> 2.5\cdot 2^n$. \epro \proofr
Рассмотрим ранги ретрактов унитрейда $U$ по некоторой координате. По
следствию \ref{corrang} суммы рангов двух ретрактов не меньше $4$.

1) Если унитрейд $U$ имеет пустой ретракт, то по утверждению
\ref{probit21} унитрейд $U$ он разложим.

2) Если унитрейд $U$ имеет  два ретракта ранга не меньше $3$ и ещё
один ранга не меньше $1$,
 то
$|U|> 2\cdot 2^n +2^{n-1}$ (см. утверждение \ref{probit101} и
следствие \ref{probit1020}).

3) Пусть все три ретракта унитрейда $U$ по координатам $x_1$ и $x_2$
имеют ранг $2$. Тогда $f=x^a\oplus x^bx_1 \oplus x^cx_2 \oplus
x^ex_1x_2$. Подставляя $x_1=1$ и $x_2=1$, получаем, что для мономов
$x^a$, $x^b$, $x^c$, $x^e$, не содержащих переменные $x_1$, $x_2$, возможны
варианты: $a\oplus b= c\oplus e$ или  $a\oplus c=b\oplus e$ или
$d(a,c)=d(a,b)=d(c,e)=d(b,e)=1$. Во всех трёх случаях имеем
$U=U[f]=U[g]U[h]$, где переменные функций $g$ и $h$ не пересекаются.
  \proofend

Для набора  векторов $W\subset Q_3^n$ определим $r(W)$ как число
позиций, в которых все наборы из $W$ совпадают. Будем считать, что
$r(W)=0$, если наборы различаются во всех позициях и в каждой
позиции принимают два значения, если же при этом найдётся позиция, в
которой присутствуют все различные буквы, то определим
$r(W)=-\infty$, т.\,е. $2^{r(W)}=0$.

\bpro\label{bitrade10} Пусть $f(x)=\bigoplus\limits_{v\in V} x^v$,
где $V\subset Q_3^n$. Тогда
$|U[f]|=\sum\limits_{t=1}^{|V|}(-2)^{t-1}\sum\limits_{W\subset V,
|W|=t}2^{r(W)}$. \epro \proofr Известна формула подобная формуле
включения--исключения:\\ $|\supp
(\chi_{A_1}\oplus\cdots\oplus\chi_{A_s})|=\sum\limits_i
|A_i|-2\sum\limits_{i\neq j} |A_i\cap A_j|+2^2\sum\limits_{i\neq j
\neq k} |A_i\cap A_j\cap A_k|-\cdots+(-2)^{s-1}|A_1\cap\cdots\cap
A_s|$.
 Как было отмечено в \S\,\ref{s:bool} справедливы равенства\\
$(U[x^{v^1}\oplus x^{v^2}\oplus\cdots\oplus
x^{v^s}])=\supp(U[x^{v^1}]\oplus\cdots\oplus U[x^{v^s}])$ и
$U[x^{v}]=\chi_{_{\{0,\pm 1\}_v}}$.

Нетрудно видеть, что $|\{0,\pm 1\}_{v^1} \cap\dots \cap \{0,\pm
1\}_{v^s}|=2^{r(\{v^1,\dots,v^s\})}$. Тогда требуемая формула
следует из подстановки этого равенства в первую формулу.
 \proofend

Пусть $v^i\in \{0,\pm 1\}^n$. Рассмотрим таблицу $\{v^i_j\}$ размера
$3\times n$, строчками которой являются векторы $v^i\in V$, $|V|=3$.
Пусть таблица содержит  $k_1(V)$ столбцов вида $acc$, $k_2(V)$ ---
$cac$, $k_3(V)$
--- $cca$ и $k_4(V)$ столбцов, состоящих из всех различных символов.
Тогда из утверждения \ref{bitrade10} имеем

\bpro Пусть $V\subset Q_3^n$, $\rank(U[f^V])=3$ и все три монома не
совпадают ни в какой координате. Тогда $|U[x^{v^1}\oplus
x^{v^2}\oplus x^{v^3}]|= 3\cdot 2^n -
2(2^{k_1(V)}+2^{k_2(V)}+2^{k_3(V)})+4\delta(k_4(V))$, где
$\delta(k)=0$, если $k>0$, и $\delta(k)=1$, если $k=0$. \epro

Рассмотрим все возможные битрейды ранга 3 мощности до $2.5\cdot2^n$
включительно.

 Заметим, что если $d(u,v)=1$, то $x^u\oplus x^v=
x^{w}$ для некоторого $w$. Поэтому достаточно рассматривать случай,
когда $k_i\leq n-2$, где $i=1,2,3$ и $n=k_1+k_2+k_3+k_4$.

1. Пусть $k_1=\max_{i=1,2,3}{k_i}=n-2$. Тогда возможны следующие
наборы

1.1. $(k_1,k_2,k_3,k_4)= (n-2,2,0,0)$. Битрейд по утверждению~\ref{bitrade91}(a). Мощность битрейда $3\cdot 2^n - 2(2^{n-2}
+4+1)+4=2.5\cdot2^n- 6$. При $n=3$ ранг битрейда равняется $2$.

1.2. $(k_1,k_2,k_3,k_4)= (n-2,1,1,0)$. Не битрейд по утверждению~\ref{bitrade91}(b).

1.3. $(k_1,k_2,k_3,k_4)= (n-2,0,0,2)$. Не битрейд по следствию~\ref{bitrade98}.

1.4. $(k_1,k_2,k_3,k_4)= (n-2,1,0,1)$. Битрейд по утверждению
\ref{bitrade9}. Мощность битрейда $3\cdot 2^n - 2(2^{n-2}+2+1)=
2.5\cdot2^n- 6$. При $n=3$ ранг битрейда равняется 2.

2. Пусть $k_1=\max_{i=1,2,3}{k_i}=n-3$. Тогда возможны следующие
наборы.

2.1. $(k_1,k_2,k_3,k_4)= (n-3,3,0,0)$. Битрейд по утверждению~\ref{bitrade91}(a). Мощность битрейда $3\cdot 2^n - 2(2^{n-3}
+8+1)+4=2.5\cdot2^n+ 2^{n-2}- 14$. При $n=4$ битрейд имеет ранг $2$,
при $n=5$ совпадает со случаем 1.1,  при $n>5$ мощность больше
$2.5\cdot2^n$.

2.2. $(k_1,k_2,k_3,k_4)= (n-3,2,1,0)$.  Не битрейд по утверждению~\ref{bitrade91}(b).

2.3. $(k_1,k_2,k_3,k_4)= (n-3,2,0,1)$.  Битрейд по утверждению~\ref{bitrade9}.  Мощность битрейда $3\cdot 2^n - 2(2^{n-3}+4+1)=
2.5\cdot2^n+ 2^{n-2}-10$. При $n=4$ совпадает со случаем 1.4, при
$n=5$ имеем мощность $2.5\cdot2^n-2$, при $n>5$ мощность больше
$2.5\cdot2^n$.

2.4. $(k_1,k_2,k_3,k_4)= (n-3,1,1,1)$.  Битрейд по утверждению~\ref{bitrade9}.   Мощность битрейда $3\cdot 2^n -
2(2^{n-3}+2+2)=2.5\cdot2^n+ 2^{n-2}-8$.  При $n=4$ имеем мощность
$2.5\cdot2^n-4$, при $n=5$ имеем мощность $2.5\cdot2^n$, при $n>5$
мощность больше $2.5\cdot2^n$.

2.5. $(k_1,k_2,k_3,k_4)= (n-3,0,1,2)$. Не битрейд по следствию~\ref{bitrade98}.

2.6. $(k_1,k_2,k_3,k_4)= (n-3,0,0,3)$. Битрейд по утверждению~\ref{bitrade9}.   Мощность битрейда $3\cdot 2^n -
2(2^{n-3}+1+1)=2.5\cdot2^n+ 2^{n-2}-4$. При $n=3$ имеем мощность
$2.5\cdot2^n-2$, при $n=4$ имеем мощность $2.5\cdot2^n$, при $n>4$
мощность больше $2.5\cdot2^n$.

 Если $\max_{i=1,2,3}{k_i}<n-3$, то мощность унитрейда больше
$2.5\cdot2^n$.

Теперь рассмотрим разложимые битрейды (см. утверждение~\ref{probit20}).

3.  Если оба сомножителя в декартовом произведении не являются
булевыми гиперкубами, то по утверждению~\ref{probit101} их мощности
могут принимать значение $\frac322^n$, $\frac742^n$ и т. д.
Поскольку $\frac32\frac74>\frac52$, мощность не более $2.5$ от
минимальной, а именно $2.25\cdot 2^n$, $n\geq 4$, имеют только
произведения битрейдов 1-го типа.

Суммируя  проведённый выше перебор и учитывая возможность декартова
умножения на булев гиперкуб (см. утверждение~\ref{probit101}),
делаем вывод, что справедлива

\begin{theorem}[мощности малых битрейдов]\label{th:2.5}
 В гиперкубе $Q^{n+m}_3$ при любом $m\geq 0$ имеются только следующие битрейды мощности
 более $2^{n+m+1}$ и  не более $5\cdot2^{n+m-1}$:\\
i) $ 2^m(2.5\cdot2^{n}- 6)$ при любых $n\geq 4$,
(1.1) и (1.4);\\
ii)  $2^m(2.5\cdot2^5-2)$  при  $n= 5$, (2.3);\\
iii)  $2^m(2.5\cdot2^3-2)$  при  $n= 3$, (2.4), (2.6) и (3);\\
iv) $2^m(2.5\cdot2^4)$  при  $n= 4$,    (2.4), (2.6).
\end{theorem}

Рассмотрим унитрейд $U$ мощности не более $2.5\cdot 2^n$ в $Q^n_k$,
$k>3$. Если по одному из направлений он  пересекается с не менее чем
с четырьмя гиперплоскостями, то пересечение с каждой гиперплоскостью
имеет мощность  $2^{n-1}$  или $3\cdot 2^{n-2}$ (см. утверждение
\ref{probit101}). Тогда мощность унитрейда $U$ может равняться
 $2\cdot 2^n$,  $2.25\cdot 2^n$ или $2.5\cdot 2^n$. Как показано выше, битрейды
 мощностей $2\cdot 2^n$ и $2.5\cdot 2^n$ имеются в троичных
 гиперкубах. Помимо них в гиперкубах $Q^n_k$ при $k>3$ имеются
 битрейды той же мощности, состоящие из двух непересекающихся
 компонент; а также битрейды вида $U=U'\times \{0,1\}^{n-2}$, где
 $U'\subset Q_k^2$ --- цикл длины 8 или 10. В $Q_4^3$ нетрудно
 построить битрейд мощности $2.25\cdot 2^3=18$. Из утверждения \ref{probit21}
 следует, что в гиперкубах $Q_4^n$, $n\geq 3$, имеются битрейды
 мощности $9\cdot 2^{n-2}$.


\section{Число битрейдов}\label{s:N}

\subsection{Нижняя оценка числа битрейдов}\label{ss:lb}

Вначале выясним, какова максимальная мощность подмножества гиперкуба
$Q^n_k$, если все попарные расстояния между  его элементами
нечётные. Начнем с рассуждений, касающихся набора вершин в евклидовом пространстве.

Пусть $\{v_1,\dots,v_n\}\subset \mathbb{R}^m$ и квадраты попарных
евклидовых расстояний между  векторами $v_i$ и $v_j$ равны
$d_{ij}=\|v_i-v_j\|_2^2$. Определителем Кэли--Менгера называется
$$\det\left(\begin{array}{ccccc}
0& 1&     1&       \dots&1  \\
1& d_{11}& d_{12}&\dots& d_{1n} \\
1& d_{21} &d_{22}&\dots& d_{2n} \\
\vdots& \vdots& \vdots & \ddots & \vdots\\
1& d_{n1} & d_{n2} & \dots & d_{nn}
\end{array} \right) = (-1)^{n+1}2^n(n!)^2(\mathrm{Vol}_{n-1})^2,$$
где $\mathrm{Vol}_{n-1}$ --- $(n-1)$-мерный объём выпуклой оболочки множества
$\{v_1,\dots,v_n\}$. 
Доказательство этой формулы объёма через детерминант 
можно найти, например, в монографиях \cite[\S\,40]{Blumenthal}, \cite[\S\,4.7]{Pak2008}. Нам важно, что при $n-1>m$  определитель равен нулю, поскольку $\mathrm{Vol}_{n-1}=0$. 
Из свойств определителя
 можно вывести следующую известную лемму.

\blemma Пусть $A\subset \mathbb{R}^m$, все попарные квадраты
евклидовых расстояний между точками множества $A$ целые нечётные и
$|A|=m+2$. Тогда $(m+2)\equiv0\mod 4$. \elemma 

\proofr  Пусть
$n=|A|=m+2$. Проведём несколько операций сложения строк и столбцов,
не меняющих определитель. Отнимем $1$-ю строку матрицы от
остальных. Получим равенство
$$\det\left(\begin{array}{ccccc}
0& 1&     1&       \cdots&1  \\
1& -1& c_{12}&\cdots& c_{1n} \\
1& c_{21} &-1&\cdots& c_{2n} \\
\vdots& \vdots& \vdots & \ddots & \vdots\\
1& c_{n1} & c_{n2} & \dots & -1
\end{array} \right) = 0,$$
где числа $c_{ij}=d_{ij}-1$ чётные.  Теперь прибавим сумму столбцов
от $2$-го до $(n+1)$-го к $1$-му столбцу, а затем  сумму строк от $2$-й до
$(n+1)$-й к $1$-й строке. Получим равенство
$$\det\left(\begin{array}{ccccc}
b& a_1&     a_2&       \dots&a_n  \\
a_1& -1& c_{12}&\dots& c_{1n} \\
a_2& c_{21} &-1&\dots& c_{2n} \\
\vdots& \vdots& \vdots & \ddots & \vdots\\
a_n& c_{n1} & c_{n2} & \dots & -1
\end{array} \right) = 0,$$ где
$a_i=\sum_j c_{ij}=\sum_j c_{ji}$ и $b=n+\sum_i
a_i=n+2\sum_{i<j}c_{ij}$. Любая диагональ матрицы,
кроме главной, содержит как минимум два чётных числа, поэтому произведение
элементов диагонали кратно $4$.
Следовательно,
произведение элементов главной диагонали также должно делиться на $4$.
Следовательно, $n\equiv b\equiv 0\bmod 4$. \proofend

\bcorol\label{bitcorol56} 
Пусть $A\subset \mathbb{R}^m$ и все
попарные квадраты евклидовых расстояний между точками множества $A$
целые нечётные. 
Тогда $|A|\leq m+2$. \ecorol 
\proofr Докажем от
противного. Пусть имеется такой набор из $m+3$ точек в
$\mathbb{R}^m$. Тогда  он также содержится в $\mathbb{R}^{m+1}$ и
удовлетворяет условию леммы, т.\,е. 
$m+3\equiv 0\bmod 4$. Кроме того, его
поднабор из $m+2$ точек в $\mathbb{R}^{m}$ также удовлетворяет
условию леммы и $m+2\equiv 0\bmod 4$. \proofend

\bcorol\label{c:(q-1)m+2} 
(a) Пусть $A\subset Q_k^m$ и все попарные  расстояния Хэмминга
между точками множества $A$  нечётные. 
Тогда $|A|\leq (q-1)m+2$.

(b) Пусть $A\subset Q_k^m$ и для любых трех точек из $A$ сумма попарных расстояний нечетна. 
Тогда $|A|\leq (q-1)m+3$.
\ecorol 
\proofr
(a) Закодируем элементы $Q_k$ наборами действительных чисел длины $k-1$ с попарными евклидовыми расстояниями $1$.
При этом словам из $Q_k^m$ будут сопоставлены векторы $(q-1)m$-мерного евклидового пространства, причём расстояние Хэмминга между словами равно квадрату евклидового расстояния между соответствующими векторами.

(b) Рассмотрим произвольную точку $a$ из $A$ и обозначим через $A'$ множество точек из $A$ на нечётном расстоянии от $a$ и через $A''$  множество точек из $A\backslash\{a\}$  на чётном расстоянии от $a$. 
Легко видеть, что расстояние между $b$ и $c$ нечётно, если $b,c\in A'$ или $b,c\in A''$, и чётно, если $b\in A'$, $c\in A''$ или $c\in A'$, $b\in A''$,
Всем наборам из $A'\cup\{a\}$ припишем в конце $0$,
а наборам из $A''$  припишем в конце $1$. Получим множество
$B\subset Q_k^m \times Q_2$ с попарно нечетными расстояниями.
Применяя далее технику, аналогичную (a), получаем множество точек в евклидовом пространстве с попарно нечетными квадратами расстояний. 
Поскольку для кодирования значений последней координаты достаточно всего одной евклидовой координаты, оценка получается на $1$ больше чем в случае (a). 
\proofend

 Для случая, когда $q$ --- степень
простого числа, общеизвестно

\bpro\label{probit11} В гиперкубе $Q^m_q$, при
$m=\frac{q^t-1}{q-1}$, имеется эквидистантный код $H_t$ мощности
$(q-1)m+1=q^t$ с кодовым расстоянием $q^{t-1}$, дуальный к коду
Хэмминга. \epro

Таким образом, при нечетном $q$, степени простого числа, есть коды, на которых граница следствия~\ref{c:(q-1)m+2} почти достигается: мощность кода всего на $1$ меньше.
Это показывает, что предложенным ниже методом 
нельзя сильно улучшить оценку снизу
на число неэквивалентных битрейдов, 
найдя код большей мощности с попарно нечетными расстояниями. 
Далее нам понадобится код $H_t$ при $q=3$.

\bpro\label{probit13} Пусть $D$ --- кодовое расстояние множества
$V_i\subset Q^n_3$ и $|V_i|2^{n-D+1}<2^{n-2}$ при $i=1,2$.  Тогда из
эквивалентности унитрейдов $U[f^{V_1}]$ и $U[f^{V_2}]$ следует
эквивалентность множеств $V_1$ и $V_2$. \epro \proofr Пусть унитрейд
$U[f^{V}]$ эквивалентен унитрейду $U'$. Тогда $U'=U[f^{V'}]$, где
множество $V'$ эквивалентно множеству $V$. Однако, соответствие
неоднозначно, т.\,е. в общем случае имеются другие  множества
$W\subset Q^n_3$, для которых $U'=U[f^{W}]$.

Достаточно показать, что если $2^{n-2}> |V|2^{n-D+1}$, то множество
$V$ с кодовым расстоянием $D$ восстанавливается однозначно по
унитрейду $U[f^V]$. Рассмотрим произвольный подкуб $U[x^v]$. Имеем
$v\in V$, если и только если $|U[f^V]\cap U[x^v]|\geq
2^{n}-2^{n-2}$. Действительно, поскольку
$|U[x^w]\cap U[x^v]|= 2^{n-d(v,w)}$ для любого $w\in Q^n_3$, имеем неравенства\\
1) $|U[f^V]\cap U[x^v]|\geq 2^n-|V|2^{n-D}$ при $v\in V$;\\
 2) $|U[f^V]\cap U[x^w]|< 2^{n-1}+|V|2^{n-D+1}$ при $w\not\in V$.
\proofend

 Обозначим через $\mathrm{sp}(v)$ состав вектора $v$, например,
 $\mathrm{sp}(0,1,1,0,-1)=(2,2,1)$. Будем говорить, что состав вектора {\it уникальный} для  некоторого линейного
 пространства, если в нём нет других векторов с тем же составом.

\bpro\label{probit12} Пусть $W\subset Q^n_k$ --- линейное
подпространство и в $W$ имеется базис $B$, состоящий из векторов с
уникальным составом. Тогда в $W$ имеется не менее $2^{|W|-\dim W
-1}/|W|$ неэквивалентных подмножеств векторов.\epro \proofr
Рассмотрим  подмножества $C\subset W$, которые содержат  нулевой
вектор и базис, т.\,е. $B\subset C$ и $\overline{0}\in C$. Пусть
 $\varphi_{\pi,a}$ --- изометрия, переводящая одно такое множество $C$ в
другое $C'$, т.\,е. $\varphi_{\pi,a}(C)=\pi(C)+a=C'$. Поскольку
$\pi(\overline{0})=\overline{0}$ имеем $a\in C'\subset W$.
Рассмотрим базисный вектор $v\in C\cap C'$ с уникальным составом. Из
равенства $\mathrm{sp}(\pi(v))=\mathrm{sp}(v)$ и уникальности состава имеем
$\pi(v)=v$. Из линейности автотопии $\pi$, т.\,е. из равенства
$\pi(\alpha u+\beta w)=\alpha\pi(u)+\beta\pi(w)$, следует, что $\pi$
действует тождественно на $W$. Тогда $\varphi_{\pi,a}(u)=u+a$, где
$a\in W$. Очевидно, что число подможеств в $W$, содержащих некоторый
базис и нулевой вектор, равно $2^{|W|-\dim W -1}$, причём любой
класс эквивалентности подмножеств содержит не больше элементов, чем
$|W|$. \proofend


Рассмотрим порождающую матрицу $A$ кода $H_t$ размерности $t>1$.
Матрица $A$ содержит единичную подматрицу. Проведём следующее
преобразование порождающей матрицы $A$: добавим к ней столбцы
единичной матрицы в количестве $2^{k-1}$ копий  $k$-го столбца при
$k=2,\dots,t$.
 Линейный код, порождённый  таким способом преобразованной
матрицей, обозначим через $H'_t$.   Все векторы кода $H_t$ (за
исключением нулевого) имеют одинаковый нечётный вес, поэтому
разность между числом координат равных $1$ и $-1$ нечётная, а значит
добавление чётного числа столбцов в порождающую матрицу не может
привести к одинаковому составу у пар коллинеарных векторов из
$H'_t$. Неколлинеарные векторы из $H'_t$ имеют разный вес по
построению. Поэтому код $H'_t$ имеет базис (строки порождающей
матрицы) из векторов с уникальным составом. Расстояние между любой
парой векторов из кода $H_t$ нечётно. Поскольку в  порождающую
матрицу $A$ было добавлено чётное число копий единичных столбцов,
расстояния между любой парой векторов из кода $H'_t$ также нечётно.
 Длина
кода $H'_t$ равна $2^t-2+\frac{3^t-1}{2}$.

\btheorem[нижняя граница] \label{th:lb} Число неэквивалентных битрейдов размерности $n$ не меньше
$2^{(2/3-o(1))n}$ при $n\rightarrow\infty$. \etheorem \proofr
 При $2^t-2+\frac{3^t-1}{2}\leq n <
2^{t+1}-2+\frac{3^{t+1}-1}{2}$ рассмотрим множество $H'_t$. По
утверждению \ref{probit11} кодовое расстояние множества $H'_t$ равно
$ D=3^t$. Для достаточно больших $t$ имеем
$|H'_t|2^{n-D+1}=3^t2^{n-D+1}\leq
3^t2^{2^t-1-\frac{3^t+1}{2}}<2^{n-2}$. Из утверждения \ref{probit13}
следует, что эквивалентность двух битрейдов $U(f^V)$ и $U(f^W)$
равнозначна эквивалентность множеств $V,W\subset H'_t$. Из
утверждения \ref{probit12} следует требуемая оценка числа таких
подмножеств.\proofend

\subsection{Верхняя оценка числа битрейдов}\label{ss:ub}

  Семейство функций $\mathcal{A}_n=\{A^n_j\}$,
  $A^n_j\subseteq \{f:Q^n_k\rightarrow S\}$,
  $n\in \mathbb{N}$, будем
  называть {\it наследственным}, 
  если оно замкнуто относительно изометрий пространства $Q^n_k$ и любой ретракт
  любой функции из $\mathcal{A}_n$ лежит в $\mathcal{A}_{n-1}$.
  Множество $T\subset Q_k^m$ называется {\it тестирующим} для
  семейства функций $\mathcal{A}_m$, если для любых $f,g\in\mathcal{A}_m$ из
  $f|_T=g|_T$ следует $f=g$.
  Множество $T$ является тестирующим для семейства $\mathcal{A}_m$ тогда и
  только тогда, когда его дополнение
  не включает носитель разности никаких двух функций из $\mathcal{A}_m$,
  т.\,е. для любых $f$ и $g$ из $\mathcal{A}_m$ имеем $Q_k^m \backslash \supp (f-g)\not\subseteq T$.
  Поскольку разность двух характеристических функций некоторых комбинаторных
  конфигураций является битрейдом (в широком смысле), то поиск
  тестирующих множеств эквивалентен нахождению множеств,
  не включающих битрейдов.

\bpro\label{probit15} Пусть семейство
  $\mathcal{A}_n$, $n\in \mathbb{N}$, наследственное. Пусть $T\subset Q_k^{m}$ ---
  тестирующее множество для $\mathcal{A}_{m}$. Тогда декартово
  произведение тестирующих множеств $T^l\subset  Q_k^{lm}$ является
  тестирующим для $\mathcal{A}_{lm}$.
\epro \proofr Докажем утверждение по индукции. Пусть
$f|_{T^l}=g|_{T^l}$. Тогда по предположению индукции для любого
$v\in T$ из $f|_{T^{l-1}\times\{v\}}=g|_{T^{l-1}\times\{v\}}$
следует, что
$f|_{Q^{(l-1)m}_k\times\{v\}}=g|_{Q^{(l-1)m}_k\times\{v\}}$. Тогда
для любого $w\in Q^{(l-1)m}_k$ имеем $f|_{\{w\}\times T}=
g|_{\{w\}\times T}$. Множество  $\{w\}\times T$ является тестирующим
для ретрактов на $\{w\}\times Q^m_k$, поскольку семейство
$\mathcal{A}_n$ наследственное. Тогда $f|_{\{w\}\times Q^m_k}=
g|_{\{w\}\times Q^m_k}$ для любого $w\in Q^{(l-1)m}_k$. \proofend

Из определения тестирующего множества и утверждения \ref{probit15}
следует

 \bpro\label{probit16} Пусть $\mathcal{A}_n$
--- наследственное семейство функций и $T\subset Q_k^m$ ---
тестирующее множество. Тогда $|\mathcal{A}_{lm}|\leq |S|^{|T|^l}$.
\epro

Ниже мы не будем различать унитрейды и их характеристические
функции. Семейства битрейдов и унитрейдов являются наследственными
(см. утверждения \ref{probit00} и \ref{probit000}). Как следует из
формулы $(\ref{eqbit1})$, тестирующим множеством для семейства
троичных унитрейдов (и битрейдов) является любое подмножество в
$Q^n_3$, индуцирующее подграф, изоморфный булеву гиперкубу.
 Пусть $T$ --- тестирующее
множество для семейства  унитрейдов в $Q^n_3$. Поскольку  число
унитрейдов в $Q^n_3$ равняется $2^{2^n}$, из утверждения
\ref{probit16} следует, что  $|T|\geq 2^n$. Отметим, что для любого
тестирующего множества $T$ его дополнение $Q^n_3\setminus T$ не
включает (непустой) унитрейд и наоборот, если $Q^n_3\setminus T$
включает унитрейд, то множество $T$ не является тестирующим для
унитрейдов. Поэтому максимальная мощность подмножества в $Q^n_3$,
не включающего унитрейды, равна $3^n-2^n$. 
Для семейства битрейдов
аналогичный вопрос остаётся открытым. 
Ниже мы по-существу
доказываем, что найдётся подмножество в $Q^n_3$  мощности больше
$3^n-2^n$, не включающее симметрические разности битрейдов.

\bpro\label{probit14} Если найдётся унитрейд $U\subset Q_3^m$,
характеристическая функция которого не является суммой (по модулю 2)
двух битрейдов, то для  битрейдов в $Q_3^m$ найдётся тестирующее
множество мощности $2^m-1$. \epro \proofr Каждой вершине $v\in
Q_3^m$ поставим в соответствие переменную $x_v$.
 Рассмотрим следующую
систему булевых уравнений однозначно задающих унитрейд $U$:\\ $
x_a\oplus x_b\oplus
x_c=0$, если вершины $a,b,c \in Q_3^m$ составляют 1-мерную грань,\\
$x_v=0, v\not\in U$.\\
Выберем из уравнений 1-го типа независимую подсистему (I), а затем
из уравнений 2-го типа максимальную независимую подсистему (II),
которая независима и с уравнениями 1-го типа. Множество решений
подсистемы (I) уравнений 1-го типа имеет размерность $2^m$, а
совместная система имеет размерность 1, поскольку (см. утверждение
\ref{probit05}) ни один  унитрейд не является подмножеством другого,
т.\,е.  нули одной функции не могут быть подмножеством нулей другой
характеристической функции унитрейда. Поэтому имеется $2^m-1$
уравнений в подсистеме (II), которые задаются точками множества
$T\subset Q_3^m$, $|T|=2^m-1$.

Покажем, что $T$ есть тестирующее множество для  битрейдов в
$Q_3^m$. Пусть две характеристические функции $\chi_A$ и $\chi_B$
различных битрейдов $A$ и $B$
 совпадают на множестве $T$, тогда
$\chi_U=\chi_A\oplus\chi_B$, что противоречит условию. \proofend

Вычислительный эксперимент (см. таблица 1) показывает, что число
битрейдов в $Q_3^7$ не превосходит $2^{2^6}=\sqrt{2^{2^7}}$ ---
квадратного корня из числа унитрейдов в $Q_3^7$. Тогда пар битрейдов
в $Q_3^7$ меньше чем унитрейдов. Таким образом при $n=7$ выполнены
условия утверждения \ref{probit14}. Отсюда получаем

 \bcorol
Число битрейдов в $Q_3^{7l}$ не превосходит $2^{{\alpha}^{7l}}$, где
$\alpha=(2^7-1)^{1/7}<2$.
 \ecorol

Пусть $\beta(n)$ --- число битрейдов в $Q_3^{n}$. Тогда
$\beta(n+m)\leq (\beta(n))^{2^m}$, $m=1,\dots,6$. Следовательно,
имеется аналогичная оценка числа битрейдов при произвольных $n>7$.

 \btheorem[верхняя граница]\label{corlast}
Число битрейдов в $Q_3^{n}$ не превосходит $2^{{\alpha_1^n}}$, где
$\alpha_1<2$. \etheorem

\end{document}